# Extracting a Paradox by the Roots

N. L. Bushwick
Northeast Autism Center, Scranton, PA, USA
nlbushwick@northeast-autism.org

**Abstract:** Zeno's paradoxes are explained as being the result of inappropriate combinations of discrete and continuous mathematical systems. It is proposed that the historical background of this confusion lies in the course of development of the number system, which was originally created to model discrete elements of experience such as apples and people, and only later, by the invention of standards of measurement, was expanded to be applied to continuous entities such as water, grain and distance. This also made it possible to combine number theory and geometry. While this modification of the system was generally unproblematic, it led to subtle contradictions when applied to time, space and motion. Analysis of these contradictions and the various resolutions that have been offered for them furthers our understanding of number systems and their relationship to physical experience. It is argued that in spite of the success of discrete systems, a truly continuous mathematical system, such as that of Carathéodory, is a more appropriate model of the continuous aspects of the physical world.

**Introduction**

The number system originated to deal with discrete objects, and only later was it adapted to continuous entities. Integers were first invented to compare collections of discrete objects such as apples, sheep and people. By counting, it became possible to determine whether there were enough horses for everyone to ride on without having to stand each person next to a horse. It became possible for a shepherd to tell whether all his sheep were back from the pasture without having to remember the name of each one. First, a list of number words such as "one, two, three" and so on had to be memorized. These are inherently arbitrary signs, whose significance resides only in their being memorized in a specific order (Wittgenstein, 1958). By putting the members of discrete collections into one-to-one correspondence with these previously memorized number words, that is, by "counting", one is able to compare their sizes.

But what about continuous entities such as rope, water and flour? To compare such things, numbers were not sufficient. It was still necessary to compare them directly, as by placing two pieces of rope next to one another, two bags of flour on the pans of a balance scale, or to resort to indirect comparison, as by pouring two quantities of water alternately into the same container and seeing which rose higher. The invention of standards of measurement addressed this problem. Units of measurement turn continuous quantities into discrete ones, making it possible to count *meters* of rope, *kilograms* of flour and *liters* of water. With the help of units of measurement, everything that could be done with discrete things could be done with continuous ones as well. The application of the number system had been extended. This new application then led to the expansion of the number system itself by motivating the invention of rational numbers to measure quantities smaller than the established units.

Time and its measurement lie somewhere in between. Time is experienced as continuous, but unlike water and rope, there are natural units. Days, months and years are all derived from natural phenomena that occur at regular intervals, and therefore lend themselves to serve as a basis by which to measure other changes. These phenomena, the cycles of the sun and moon, are such a striking part of human experience that it is impossible to fail to notice them. Even the most primitive cultures therefore measure the passage of days, months and years. That is to say, there is a natural discreteness about the experience of time. Time also differs in yet a more essential way. Time is not a material. Indeed, it is not even a phenomenon, but rather an abstract concept. What is experienced is *change*, not the flow of some material called "time". There is always the *experience* of the present, the *memory* of the past, and the *expectation* of the future. It is the human mind that synthesizes them and creates the concept of time (Carr, 1912; Harré, 1996; cf. also von Wright, 1969).

But, while the invention of standards of measurement is an excellent practical solution to the problem of applying numbers to continuous entities, it may lead to the illusion that they are *in fact* collections of discrete objects. Measurement tends to conceal the essential difference between the equality of discrete entities and that of continuous ones. Discrete equality is, by definition, exact. It is arrived at by putting two discrete collections into one-to-one correspondence with each other. Equality of continuous entities may be exact, but more often it incorporates a degree of tolerance of deviation. If there are ten people and ten horses, then each person gets a horse to ride. Exactly one. There are no horses left over and no one has to walk. But if there are ten quarts of water, each person need not get exactly one. Units of measurement make continuous entities *functionally* discrete but not *actually* so. They make it possible to treat them as if they were, but do not make them so in fact. Measurement is meant to be a practical technique. It is not designed to reveal essential truths about the physical world. The illusion that continuous entities are in fact collections of discrete objects is of no practical harm. Commerce and even dosages of medicine work even if the measurements are not exact. But when it comes to philosophical arguments, this can lead to serious errors.

That illusion is not produced by the institution of measurement itself. Standards of measurement in themselves do not imply exactness. The ancient Babylonian merchant, who had learnt to measure grain as earlier he had been able to count sheep, would probably not have claimed that each measure of grain was exactly the same. Indeed, no two sheep are exactly the same either. All he was interested in was the practical technology of counting for the purpose of commerce. The fact that an animal is a natural unit while a measure of grain is an artifact was only incidental. All that mattered was that it rendered grain discrete and therefore countable. Standards of measurement, therefore, contain the seeds of future error, but the proper soil is yet required for them to bear fruit.

It was in classical Greek culture that the confusion arose, for there these two kinds of applied mathematics, discrete and continuous, gave rise to two corresponding branches of abstract mathematics, *number theory* and *geometry*. These were already being developed by the time of Pythagoras, and are most familiar to us through the works of Euclid, who compiled them several centuries later. The use of geometric configurations to symbolize the propositions of number theory may have contributed to conceptual crossing over between number theory and geometry. Classical mathematical notation was too crude and cumbersome to conveniently express the properties of number theory. It tended to obscure relationships between numbers rather than reveal them. The use of geometric configurations instead was a clever solution, but it also invited confusion. Counting objects in an array is a discrete operation. It has nothing to do with the sizes of the objects counted or the distances between them. It only involves their number. Geometry, on the other hand, deals with distance, area, and angle, which are continuous. But since the shapes of the arrays used in classical number theory also appear in geometry, there is a large body of analysis that is common to the two.

Within geometry itself another confusion developed. Even though the physical world is actually three dimensional, there are physical situations from which the concepts of one and two dimensional forms can be abstracted. Land area is two dimensional, and distance is one dimensional. The concept of a dimensionless point follows, for the intersection of two lines is a point. Abstract geometry involves the combination of forms of differing dimensions, lines together with planes, planes and lines with solids. Each is discrete relative to those of higher dimension. Points on a line and lines on a plane are discrete objects. Not so segments of a line or two-dimensional forms within a plane, which are continuous. The combination of objects of different dimensionality, such as points and lines, in a single system tends therefore to lead to contradictions.

The Greeks soon became aware of the contradictions into which the marriage of discrete and continuous systems was leading them. These were expressed by Zeno in his four paradoxes of motion: *the race course, Achilles and the tortoise, the arrow*, and *the stadium*. Though it is primarily the paradox of the arrow that is of interest to us here, the others need to be discussed as well. In the second and most famous, Achilles, the fastest of mortals, is to run a race with a tortoise. Since the tortoise is so slow, it is given a head start. Now, we all know that Achilles will soon catch up with the tortoise and overtake it, but Zeno argued that if distance could be continually divided and redivided, that would be logically impossible. By the time Achilles reaches the tortoise's starting point, the tortoise has moved on, even if only a very short distance, so Achilles has not yet caught up. He keeps running and reaches the spot where the tortoise was then, but again it has moved ahead, so he runs further. This happens over and over. However close they are, the relationship between them remains the same; only the distance changes. The number of 'runs' is endless. No matter how many Achilles has covered there will always be another. As there is never a last 'run', the race is never over and Achilles never catches up with the tortoise.

The paradox of the race course presents a similar problem. Here, a would-be runner is to run a course from start to finish. But before he can complete the whole course he must cover half of it, before that one quarter, before that one eighth and so on. Since there is no end to these divisions, there is no first stretch for him to cover. He is in even worse shape than Achilles, because not only does he not reach his goal, he never even gets started.

In the third paradox, Zeno has us consider an arrow flying through the air. At any moment the arrow occupies a certain single location, exactly equal in size to that of the arrow itself. At that moment there is no difference between the flying arrow and one that is stationary at that spot. Each occupies but one location at any given moment. At each moment, then, the arrow is not moving. But since the entire flight of the arrow is made up of the sum of these moments, the arrow is never moving. How, asked Zeno, can an arrow that is always standing still get from one place to another?

Although classical philosophy did not explicitly articulate the concepts 'continuous' and 'discrete', the extended discussion of which these paradoxes were a part clearly involved them. The classical concept most closely related was that of *divisibility*. The question was whether the physical world could be divided over and over again without end, or whether eventually some smallest indivisible unit would be reached. The first two paradoxes argue against infinite divisions, corresponding to a continuous view of the world, the third against minimal units. It is not clear what Zeno intended to prove by these arguments. Extensive discussion of the historical question can be found in Tannery (1987), Owen (1970), Papa-Grimaldi (1996), Arntzenius

(2000), Glazebrook (2001), and other sources. But however the paradoxes are interpreted, each involves a combination of discrete and continuous descriptions, and it is that combination that produces the paradox.

Each of these three paradoxes involves an infinite summation, the adding together of an endless collection. This is a discrete concept, involving the discrete operation of addition. The problem of how such a summation can ever be completed immediately arises. If the terms are to be added one by one, the sum can never be completed because there is no last term. In the arrow paradox, the problem is compounded by the introduction of another discrete concept, that of dimensionless points or moments. The first two paradoxes involve summations of successively smaller terms, each of which, however, has extension. Here, beside being infinitely many, the parts of which the line is imagined to consist are unextended points. This is problematic even without applying it to physical entities such as ropes or lines drawn on paper. How can a line be composed of dimensionless points, since a sum of zeros, no matter how many, can never add up to anything other than zero?

**Resolving the Paradoxes - A Discrete Resolution**

There are two basic approaches to resolving these paradoxes, and each has had adherents since antiquity. Essentially, one approach rejects the continuous description in favor of the discrete, and the other the discrete in favor of the continuous.

One early and brilliant discrete solution is provided by atomism (Owen, 1970). The classical atomists and their followers in later periods held that matter is composed of very small indivisible parts, so it cannot be divided and redivided forever. Eventually the smallest units of matter will be reached and it will not be divisible any more. The atomist description applies not only to matter but to time and space as well. Rather than being viewed as continua, they are seen as being composed of discrete moments and discrete locations that cannot be divided further. Time is *in fact* like the separate frames of a motion picture, and it is only our perception that makes it appear smooth. So too, space is like the pixels of a digital image. That we are unable to see the spaces between them is not just a practical limitation of our nervous systems, but the essential limitation of our existence, since we ourselves are composed of exactly those units of time and space. There might be miles between one unit of space and the next, and millennia between an instant and the one that follows it, yet it would be essentially unknowable for us, since we ourselves do not exist there. Indeed, the very idea of time and space between the units is meaningless, since nothing at all exists between them. It is but the error of an imagination that still considers existence to be continuous and has not fully grasped the concept of atomism (Arntzenius, 2000; cf. von Wright, 1969).

According to this position, the problems of the paradoxes are illusory. In the first two paradoxes, the description itself is fallacious. The length of the race course cannot in fact be divided endlessly. Eventually the length of the division is a single unit of space, and then the process must stop. So too, when the distance between Achilles and the tortoise is only one unit he will catch up on the next 'run'.

As for the arrow, motion is indeed an illusion, if by motion one means continuous change of position. The reality of motion is that the moving object occupies successive positions in space at successive moments in time. This is what we perceive as motion. Whether the number of positions of the arrow as it moves from one end of its path to the next is finite, as some atomists maintained, or infinite, the source of the paradox is a misunderstanding of the true nature of motion. It is not a smooth changing of position, as it appears to be, but a succession of separate abrupt changes.

It was against the atomist position that the fourth paradox, the paradox of the stadium, was directed. If there are three rows of men aligned with one another in a stadium, one marching in one direction, another at the same rate in the opposite direction, and the third standing still, what happens in the smallest unit of time? If, relative to one another, the moving rows have advanced by one minimal unit of distance, then relative to the stationary one they have only moved half a unit, which contradicts the assumption that the unit is minimal (Tannery, 1987; Owen, 1970; Glazebrook, 2001). To maintain the atomist position, one is forced to say that in the minimal time unit the two moving rows advance two units relative to one another, and that there is no time unit in which their relative motion is only one.

There seems, in principle, to be no way of refuting this form of atomism, since, if it is true, our existence is entirely within the discrete atomic world, so there is no way for us to see the transition between atoms. Furthermore, it cannot be refuted practically by making very small divisions, since however small a division we succeed in making, it can be argued that the atoms are smaller yet. On the other hand, unless some proof for atomism is devised, we are inclined to consider the perceived continuity of space and time to reflect reality.

**A Continuous Resolution**

The opposite approach was taken by Aristotle. In the Sixth Book of *Physics*, in the course of presenting a theory of space, time and matter, he resolves the paradoxes within an internally consistent continuous system. The first two paradoxes are answered by rejecting the need for a discrete additive process for the composition of a continuum. Even though any continuum can be divided endlessly, that is only a potential. Zeno presented it as if all the potential divisions had actually been made, and therefore needed to be recomposed, but Aristotle maintained that the motion of a moving body does not consist of separate events. Achilles' motion is a single smooth run, as is that of the runner on the race course. Since it is not really composed of separate 'runs', there is no need to add them together. Furthermore, however many divisions are made, they are smaller and smaller and their sum is always just that of the original distance. It is only the artificial division that causes motion to appear paradoxical (Glazebrook, 2001). Space, time and motion are really continuous, and it is only when they are erroneously described as being discrete that there is a problem. The fallacy in the paradoxes of Achilles and of the race course lies therefore in the description.

While the discrete approach accounts for the paradox of the arrow in the same way as it does the first two, for the continuous approach this third paradox requires a different kind of solution. Here, one is dealing with instantaneous moments rather than smaller and smaller but always measurable intervals. The problem is the opposite of that of the first two. How can the sum be other than zero? The problem of composing a time period from instantaneous moments is the same as that of composing a line from dimensionless points. Aristotle's answer is again to eliminate the discrete aspect of the description and maintain that even potentially, a line cannot be considered to be composed of points. Since two distinct points cannot touch one another, there must always be something between them. So even though there are infinite points on the line, the line is composed not of these points but of the distances between them (Arntzenius, 2000; Glazebrook, 2001).

**Comparison of the Two Approaches**

Atomist and continuous solutions both succeed for the same reason. In each of the four paradoxes, the apparent contradiction is the result of an inconsistent combination of discrete and continuous descriptions. The atomist solution is to revise the continuous part and make it discrete. The continuous solution is to make the whole thing continuous.

In neither case is there any attempt to explain time itself and how it progresses. The passage of time is simply a fact. It is the nature of our existence that we experience change, from which we create the concept of time (von Wright, 1969). It is essentially impossible to understand time itself, because we and our processes of observation and of thinking are entirely within it. Whether it 'ticks' or 'sweeps', it keeps moving. But lack of understanding of the essence of time neither adds to nor detracts from the force of the paradoxes, and however one resolves them, time itself remains a mystery.

These two approaches reflect two radically different attitudes toward knowledge. The discrete approach is based upon an idealist doctrine, such as that of Pythagoras, that the physical world is derived from numbers, and ultimately, from the number one (Tannery, 1987; Glazebrook, 2001). Beginning a priori with such a doctrine, experience is forced to conform to it. Since the number system is discrete, the world must be discrete too. Indeed, the paradoxes may have been directed against exactly this position, Zeno's argument being that any application of discrete mathematical operations to the continuous physical world is bound to produce a contradiction. According to this interpretation, the concepts of 'one' and 'many' used by Zeno and his teacher, Parmenides, were close, if not identical, to the concepts of 'continuous' and 'discrete' (Papa-Grimaldi, 1996). But the adherents of the discrete position, rather than relinquishing it, responded by further removing continuity from their description of the physical world.

The continuous position is derived from the physical experience that there exist some entities that can be divided and redivided, the resulting parts always being divisible themselves. Extrapolation from this experience implies the reality of continuity, and experience is given preference over any a priori system (Glazebrook, 2001). Aristotle saw that the source of the contradiction lay in the introduction of entities such as points and operations such as addition, derived from a discrete formal system, into the continuous system derived from physical experience. He addressed the problem by limiting the application of discrete entities and operations, in the one case by asserting that a continuum could not be composed of points, and in the other by denying the necessity of performing infinitely many discrete additions to compose a whole. He did not, however, exclude discrete aspects entirely, limiting them, rather, only as much as he considered necessary to avoid contradiction.

**Seeing the Problem as One of Representation**

The ancient philosophers considered the subject of the paradoxes to be the physical world itself, the question being whether it was endlessly divisible or not. Their arguments, however, dealt not with the physical world but with descriptions of it, their goal being to prove or disprove the consistency of those descriptions. Indeed, the same can be said of most later discussions of the paradoxes. The real question is therefore whether or not it is possible to construct a description of the physical world that is internally consistent. The paradoxes argue that it is not, because any description will contain elements that contradict one another. The resolutions are achieved by modifying the descriptions to make them internally consistent (Salmon, 1970).

Viewing the problem as one of representation is a radical departure from the traditional interpretation of its being a question about the physical world itself. There are now two separate questions, one being whether a particular description is internally consistent and the other being whether it is an appropriate description of the physical world. Is that the way the world really is? These two questions are, of course, not independent. If a description is not internally consistent it cannot be appropriate, since the world itself is not inconsistent. However, internal consistency alone does not make it so (Salmon, 1970).

This view implies, moreover, that there may be more than one valid resolution of the paradoxes. Consistency and appropriateness alone do not guarantee that a description is unique. There may be more that one consistent and appropriate description of the physical world. It is well known that there are many things that can be described in more than one way, each way being an appropriate description as well as being internally consistent. Sometimes, different descriptions address different aspects. For example, political and topographic maps represent different aspects of the same continent. In other cases various descriptions address the same aspect, but in different ways. A complex number, for example, can be described using either polar or rectangular coordinates. For the paradoxes too, therefore, there may be two or more resolutions.

Those, past and present, who imagine that by achieving consistency they have identified the unique true nature of the world have fallen into the common error of confusing description with its referent. Tannery (1987) and others after him (Salmon, 1970; Papa-Grimaldi, 1996; Glazebrook, 2001) have attributed a certain degree of this insight to Zeno himself. According to this interpretation, the purpose of Zeno's arguments was not to prove that motion was impossible, nor that the world was neither continuous nor atomic, but rather that mathematics was incapable of providing a consistent description of physical reality. What Zeno was doing, then, was demonstrating a limitation of mathematics. He was not addressing the nature of the physical world at all. Even this very generous interpretation of Zeno falls short of the complete insight, the problem being seen only as one of mathematics rather than of description in general. However, if indeed this was Zeno's intention, he came about as close to identifying the problem as one lacking the concepts of continuous and discrete systems could have. If so, however, he was not understood by subsequent philosophers, and the insight latent in it remained unnoticed for over two millennia (Glazebrook, 2001).

**The Ascent of Discrete Mathematics**

Symbols are tools of the mind, increasing its power to manipulate and analyze concepts, but in doing so they bias it to certain attitudes and modes of thinking. Since the introduction of Arabic numerals in the Middle Ages, there has been a continual advance in mathematical notation, which has facilitated the advance of mathematics itself. The decimal system made it possible to represent arbitrarily large and arbitrarily small quantities and to perform arithmetic operations conveniently, which enabled mathematicians to think of numbers without resorting to geometric representations. But since these new and powerful symbols were discrete, discrete analysis came to increasingly dominate mathematical thought. Most of the currently accepted solutions of the paradoxes are therefore discrete ones (Russell, 1993; Zangari, 1994; Papa-Grimaldi, 1996; A notable exception was Bergson, 1927, 1992; Carr, 1912; cf. also von Wright, 1969). Furthermore, the more powerful the discrete system became, the more it seemed to be the unique correct description of the physical world, and the greater the effort it inspired to model the world according to it.

During the nineteenth century, a new definition of continuity itself was developed. Largely through the work of Bernard Bolzano (1950), the concept of an extended entity was replaced by that of a dense set of points. This was particularly welcome in the development of calculus. The intuitive notion of continuity that calculus had relied on until then was problematic because it was not compatible with the discrete concept of a function, which takes discrete values and yields discrete values as its results. By reconstruing continuous entities such as space and time as dense sets of discrete elements, the concept of a continuum was converted from that of an extended entity that could be repeatedly redivided but never decomposed into elementary units to that of a collection of points, thus making it in effect a discrete concept. By this step, together with advances in understanding of convergence and of orders of infinity, it became possible to resolve the paradoxes without having to resort to minimal atomic units as classical atomism had.

Several objections were thereby avoided. Since irrational numbers are an uncountable infinity, their sum is defined neither by counting nor by a limit of counting, so mathematical induction does not apply to them. It can therefore not be argued that the sum of any number of zeros is zero. Nor are current discrete approaches subject to the Aristotelian objection that, since no two points can be adjacent, some of the line has been left out. That objection can only be made if such a missing place on the line can somehow be identified, since it is not meaningful to speak of a place on the line unless it can be indicated in some way. That identification can certainly not be done by physically pointing to it, since a line is a formal, not a physical, concept. It must therefore be done by a number, there being no other kind of description. So the line consists exactly of its points, corresponding to the real numbers, rational and irrational, and nothing has been left out. Defining a line interval as an uncountably infinite collection of points is therefore internally consistent, and can serve both as a rigorous basis for calculus and an approach to resolving the paradoxes.

The success of discrete mathematics left little motivation for the construction of a continuous system. However, in the twentieth century Constantin Carathéodory constructed a system that, without resorting to dimensionless points, could serve as an alternate basis for calculus and physics. It is internally consistent and not subject to the objections of the paradoxes. The elements of this system are extended geometric entities, which are not defined in terms of end points or boundary lines (Carathéodory, 1963; Arntzenius, 2000). In describing a line, for example, segments themselves are elementary objects. While the word "interval", as used in the traditional systems, refers to the region between two values, and must therefore be defined in terms of end points, the word "segment" is not defined in that way. Points are therefore not required, and need not be included in the system. The concept of continuity in such a system is that of an extended entity, not of a dense collection of points.

**Why Prefer a Continuous Description?**

Given, then, that there are now two kinds of mathematical models, both of which are internally consistent and both of which can be used to describe physical phenomena, which is to be preferred? At first it would seem to be the discrete ones. Since some aspects of the world, such as apples and sheep, are undeniably discrete, no system entirely lacking in discrete aspects can be sufficient to represent it. And since those aspects of the world that seem to us to be continuous can be successfully described discretely, it is possible, using a discrete system, to apply a single mathematical foundation to everything. A continuous description, on the other hand, must necessarily be a hybrid. Whatever system it uses to represent the continuous aspects, it must include a number system of some sort to represent the discrete ones. Besides being more complicated, such a combination tends to lead to inconsistency, as Zeno demonstrated. It would seem, therefore, that continuous descriptions should be avoided.

But are discrete models really appropriate descriptions of the physical world? Our experience of the world is always three dimensional and always subject to constantly changing time. It contains neither points nor durationless instants. Time is always experienced as moving and space as having extension. Points in space and instants in time are cultural illusions, engendered by discrete descriptions. They are not natural human mistakes. The stone-age artist who painted a herd of deer running across the wall of a cave would not have said that they were standing still. During the past half century this illusion has been intensified by the proliferation of digital computers and the use of digital measuring devices in place of analogue ones. The world is viewed increasingly through digital representations rather than directly, and our mode of thought and expression becomes increasingly discrete.

The isolation of separate dimensions is a theoretical concept abstracted from that experience. This act of abstraction is similar to other abstractions that we make, and that pervade our thought and language. From a triangular piece of wood we abstract the material "wood" and the shape "triangle", but never in our experience is there a triangle that is not made of some material nor wood that does not have some shape. That is not to say that shapes, materials and dimensions are not real. They are real properties of objects and serve as useful descriptions, but they are not objects themselves. They cannot exist except in certain combinations, such as shape together with material. So too, the three dimensions and time can only exist together with one another. Redefining continuity as density of points does not constitute a deeper understanding of the physical reality. What has really been accomplished is not a description of continuity, but a discrete simulation of it.

To build an argument upon the state of an arrow at a particular moment in its flight is like basing a proof upon the description of a shapeless piece of wood. Just as it is useful to speak of a material or a shape, so too it is useful to speak of locations of a moving arrow. But while we can measure a distance and describe it, the arrow never exists at that point alone. In a durationless instant there is no arrow, no world, no existence (cf. von Wright, 1969). We might take a photograph of the arrow at that point, but the photograph is not an arrow, it is only a picture of an arrow. It is a kind of description of an arrow, somewhat like a verbal description.

The concept of instantaneous velocity is useful in analyzing the physical world just as standards of measurement are useful in commerce. The fact that a calculus based upon a discrete foundation provides us with a way of defining and calculating instantaneous velocity does not mean that such a physical reality actually exists. Just as in language it is possible to make up words for imaginary objects such as dragons and unicorns, so in mathematics it is possible to create models to which nothing in the physical world corresponds (Glazebrook, 2001). The usefulness of this definition in physics is comparable to the usefulness of standards of measurement in commerce. Neither need be seen as a description of a true essence.

The popularity of discrete systems is part of the tradition that evolved in the course of the development of mathematics and physics. From the beginning, classical geometry involved interrelationships of elements of different dimensions. Such interrelationships being essential and fundamental to it, not merely incidental, their implications could not be escaped, leading ultimately to an essentially discrete system. But since the physical world involves only entities of three dimensions, such a system cannot be an accurate description of it. Though for one raised in the Western scientific tradition the discrete approach may seem natural and appropriate, when considered without this bias, a continuous system is a much better description of the physical world.

Zangari, M. 1994 "Zeno, Zero, and Indeterminate Forms: Instants in the Logic of Motion", *Australasian Journal of Philosophy* 72(2), pp. 187-204